\documentclass{amsart}

\usepackage{graphicx}
\usepackage{amsthm}
\usepackage{amsmath}
\usepackage{amsfonts}
\usepackage{amssymb}

\newtheorem{theorem}{Theorem}

\newtheorem*{corollary}{Corollary}

\newtheorem*{definition}{Definition}
\newtheorem{example}{Example}

\newtheorem*{lemma}{Lemma}

\newtheorem*{problem}{Problem}
\newtheorem{proposition}{Proposition}
\newtheorem*{remark}{Remark}

\begin{document}

\title[Products by generators]{On the product by generators of characteristically nilpotent Lie S-algebras\footnote{Research partially supported by the D.G.I.C.Y.T project PB98-0758}}
\author[J.M. Ancochea]{Jos\'{e} Mar\'{\i}a Ancochea Berm\'{u}dez\footnote{%
corresponding author: e-mail: Jose\_Ancochea@mat.ucm.es}} 
\author[R. Campoamor]{Rutwig Campoamor Stursberg\\Departamento de Geometr\'{\i}a y Topolog\'{\i}a\\
Fac. CC. Matem\'{a}ticas Univ. Complutense\\
28040 Madrid ( Spain )}
\date{}
\maketitle
\begin{abstract}
We introduce the product by generators of complex nilpotent Lie algebras, which is a commutative product obtained from a central extension of the direct sum of Lie algebras. We show that the product preserves also the characteristic nilpotence provided that the multiplied algebras are $S$-algebras. In particular, this shows the existence of nonsplit characteristically nilpotent Lie algebras $\frak{h}$ such that the quotient $\frac{\dim \frak{h}-\dim Z\left(  \frak{h}\right)}{\dim Z\left(\frak{h}\right)}$ is as small as wanted.
\end{abstract}

\bigskip
A nilpotent Lie algebra $\frak{g}$ is called characteristically nilpotent (
short CNLA) if there exists an integer $m$ such that $\frak{g}^{[m]}=0$, where
\[
\frak{g}^{[1]}=Der\left(  \frak{g}\right)  \left(  \frak{g}\right)  =\left\{
X\in\frak{g}\;|\;X=f\left(  Y\right)  ,\;f\in Der\left(  \frak{g}\right)
,\;Y\in\frak{g}\right\}
\]
and
\[
\frak{g}^{[k]}=Der\left(  \frak{g}\right)  \left(  \frak{g}^{[k-1]}\right)
,\;k>1
\]
In [8] it is shown that this condition is equivalent to the nilpotence of the Lie algebra of derivations. CNLAs were introduced in 1957 by Dixmier and Lister [4], and have
constituted an interesting object of study since then. It is well known that,
given two CNLAs $\frak{g}_{1}$ and $\frak{g}_{2}$, their direct sum
$\frak{g}_{1}\oplus\frak{g}_{2}$ is also a CNLA. This allows to construct such
algebras in arbitrary dimension. However, the disadvantage is the splitness of
the resulting algebra. In this paper we introduce a commutative operation
$\underline{\times}$ on the variety of complex nilpotent Lie algebras, which is called
product by generators. This product is a natural generalization of the direct sum, and it has the remarkable property of preserving nonsplitness. Moreover, for the subclass of nilpotent S-Lie algebras the product by generators of two CNLAs is also characteristically nilpotent ( and nonsplit), generalizing the result obtained by Leger and T\^{o}g\^{o} in [8] for the direct sum. 
In particular, we can construct nonsplit CNLAs for which the proportion between the elements in the derived algebra which are not central and the center is as small as wanted.\newline We convene that whenever we speak about nilpotent Lie algebras, we mean a finite dimensional complex nonabelian nilpotent Lie algebra.   

\bigskip
Let $\frak{g}_{1}$ and $\frak{g}_{2}$ be nilpotent Lie algebras of dimensions
$m_{1},m_{2}$, respectively. Suppose that $\left\{  X_{1},..,X_{m_{1}%
}\right\}  $ is a basis for $\frak{g}_{1}$ such that $\left\{  X_{1}%
,..,X_{n_{1}}\right\}  $ generates the algebra, and $\left\{  X_{1}^{\prime
},..,X_{m_{2}}^{\prime}\right\}  $ a basis of $\frak{g}_{2}$ with generators
$\left\{  X_{1}^{^{\prime}},..,X_{n_{2}}^{\prime}\right\}  $. Consider the
direct sum $\frak{g}_{1}\oplus\frak{g}_{2}$ and the 2-cochain $\varphi\in
C^{2}\left(  \frak{g}_{1}\oplus\frak{g}_{2},\mathbb{C}^{n_{1}n_{2}}\right)  $
given by%
\[
\varphi\left(  X_{i},X_{j}^{\prime}\right)  =\left\{
\begin{array}
[c]{cc}%
e_{\left(  i-1\right)  n_{2}+j} & 1\leq i\leq n_{1},1\leq j\leq n_{2}\\
0 & \text{if \ }i>n_{1}\text{ or \ }j>n_{2}%
\end{array}
\right.
\]
where $\left\{  e_{1},..,e_{n_{1}n_{2}}\right\}  $ is the canonical basis of
$\mathbb{C}^{n_{1}n_{2}}$.

The proof of the next assertion follows at once : 

\begin{lemma}
For $n_{1},n_{2}\geq2$ we have $\varphi\in Z^{2}\left(  \frak{g}_{1}%
\oplus\frak{g}_{2},\mathbb{C}^{n_{1}n_{2}}\right)  $.
\end{lemma}

It is well that the central extensions of a nilpotent Lie algebra $\frak{g}$ are defined by cocycles for the trivial representation [5]. This justifies the following :

\begin{definition}
The product by generators $\frak{g}_{1}\underline{\times}\frak{g}_{2}$ of
$\frak{g}_{1}$ and $\frak{g}_{2}$ is the central extension of $\frak{g}%
_{1}\oplus\frak{g}_{1}$ defined by the cocycle class of $\varphi$.
\end{definition}

\begin{proposition}
If $\frak{g}_{i}$ is nonsplit for $i=1,2$, then the product by generators
$\frak{g}_{1}\underline{\times}\frak{g}_{2}$ is a $\left(  m_{1}+m_{2}%
+n_{1}n_{2}\right)  $-dimensional nonsplit Lie algebra of nilindex
$\max\left\{  n\left(  \frak{g}_{1}\right)  ,n\left(  \frak{g}_{2}\right)
\right\}  $.
\end{proposition}

\begin{remark}
We observe in particular that \textbf{any} product by generators of two Lie algebras ( abelian or not ) is nonsplit.
\end{remark}

\begin{proposition}
For nilpotent Lie algebras $\frak{g}_{i}\;\left(  1\leq i\leq t\right)  $ we have

\begin{enumerate}
\item $\frak{g}_{1}\underline{\times}\frak{g}_{2}=\frak{g}_{2}\underline
{\times}\frak{g}_{1}$

\item $\left(  \frak{g}_{1}\underline{\times}\frak{g}_{2}\underline{\times
}\frak{g}_{3}\underline{\times}...\frak{g}_{1}\underline{\times}\frak{g}%
_{t-1}\right)  \underline{\times}\frak{g}_{t}\simeq\frak{g}_{1}\underline
{\times}\left(  \frak{g}_{2}\underline{\times}\frak{g}_{3}\underline{\times
}\frak{g}_{4}\underline{\times}...\frak{g}_{t-1}\underline{\times}\frak{g}%
_{t}\right)  $
\end{enumerate}
\end{proposition}

\begin{definition}
A nilpotent Lie algebra $\frak{g}$ is called a product by generators if there
exist two subalgebras $\frak{h}_{1}$ and $\frak{h}_{2}$ such that
$\frak{h}_{1}\underline{\times}\frak{h}_{2}$.
\end{definition}

As usual, the center of a nilpotent Lie algebra $\frak{g}$ is denoted by $Z\left(
\frak{g}\right)  $, while the ideals of the central descending sequence are
denoted by $C^{i}\frak{g}$, where $C^{1}\frak{g=[g},\frak{g]}$. 

\begin{proposition}
If $\frak{g}$ is a product by generators, then $\dim Z\left(  \frak{g}\right)
\geq6$ and $\dim\frac{\frak{g}}{C^{1}\frak{g}}\geq4$. In particular, no
filiform Lie algebra can be a product by generators.
\end{proposition}

\begin{proof}
Let $\frak{g}$ be a such product and $\frak{g}_{1},\frak{g}_{2}$ the ( nonabelian )
subalgebras which define it. As any nilpotent Lie algebra has at least two
generators, and the product by generators combines precisely these vectors,
the product $\frak{g}_{1}\underline{\times}\frak{g}_{2}$ has at least four
generators. Now the adjoined vectors are central, and as $\dim Z\left(
\frak{g}_{i}\right)  \geq1$ for $i=1,2$, the assertion follows. Now, if
$\frak{g}$ were filiform, we would have $\dim\frac{\frak{g}}{C^{1}\frak{g}}=2.$
\end{proof}

\begin{corollary}
If $\frak{g}$ is a product by generators, then $\dim\left(  \frak{g}\right)
\geq10$.
\end{corollary}

\begin{remark}
The preceding dimension is the lowest possible, as we have excluded the one dimensional Lie algebra $\mathbb{C}$ from our analysis.
\end{remark}

\begin{corollary}
If $\frak{g}$ is a product by generators of characteristically nilpotent algebras,
then $\dim\left(  \frak{g}\right)  \geq18.$
\end{corollary}

\begin{proof}
The lowest dimension for which CNLA exist is seven [6]. Thus, as these
algebras have at least two generators and a one dimensional center, the
minimal dimension for a product $\frak{g}_{1}\underline{\times}\frak{g}_{2}$
is 18.
\end{proof}

Let $\frak{g}_{1}$ and $\frak{g}_{2}$ be nilpotent Lie algebras and consider
the product by generators $\frak{g}_{1}\underline{\times}\frak{g}_{2}$. It is
clear that as vector space we have $\frak{g}_{1}\underline{\times}\frak{g}%
_{2}=\frak{g}_{1}\oplus\frak{g}_{2}\oplus\frak{g}_{3}$, where $\frak{g}%
_{3}=\left\langle Z_{1},..,Z_{n_{1}n_{2}}\right\rangle $. Let $p_{i}%
:\frak{g}_{1}\underline{\times}\frak{g}_{2}\rightarrow\frak{g}_{i}$ \ be the
projection of $\frak{g}_{1}\underline{\times}\frak{g}_{2}$ over $\frak{g}_{i}$
for $i=1,2,3$. For an arbitrary derivation $D\in Der\left(  \frak{g}%
_{1}\underline{\times}\frak{g}_{2}\right)  $ write $D_{ij}=p_{i}Dp_{j}$. It is
obvious that $D=\sum_{i,j}D_{ij}$. In particular
\begin{align*}
D\left(\frak{g}_{1}\right) &  =D_{11}\left(\frak{g}_{1}\right)+D_{21}\left(\frak{g}_{1}\right)+D_{31}\left(\frak{g}_{1}\right)\\
D\left(\frak{g}_{2}\right) &  =D_{12}\left(\frak{g}_{2}\right)+D_{22}\left(\frak{g}_{2}\right)+D_{32}\left(\frak{g}_{2}\right)\\
D\left(\frak{g}_{3}\right) &  =D_{13}\left(\frak{g}_{3}\right)+D_{23}\left(\frak{g}_{3}\right)+D_{33}\left(\frak{g}_{3}\right)%
\end{align*}

\begin{proposition}
The following relations hold :

\begin{enumerate}
\item $D_{ii}$ is a derivation of $\frak{g}_{i}$ for $i=1,2$

\item $D_{21}\left(  C^{1}\frak{g}_{1}\right)  =D_{12}\left(  C^{1}%
\frak{g}_{2}\right)  =0$

\item $D_{i3}\left(  \frak{g}_{3}\right)  \subset Z\left(  \frak{g}%
_{i}\right)  \;$ for $i=1,2$

\item $D_{12}\left(  \frak{g}_{2}\right)  \subset T_{1},\;D_{21}\left(
\frak{g}_{1}\right)  \subset T_{2}$, where%
\[
T_{i}=\left\{  Y\in\frak{g}_{i}\;|\;\left[  Y,\frak{g}_{i}\right]  \subset
Z\left(  \frak{g}_{i}\right)  \right\}  ,\;i=1,2
\]
\end{enumerate}
\end{proposition}

\begin{proof}
Let $X,X^{\prime}$ be arbitrary elements of $\frak{g}_{1}$. Then
\begin{align*}
D\left[  X,X^{\prime}\right]   &  =D_{11}\left[  X,X^{\prime}\right]
+D_{21}\left[  X,X^{\prime}\right]  +D_{31}\left[  X,X^{\prime}\right]  \\
&  =\left[  D_{11}X,X^{\prime}\right]  +\left[  D_{21}X,X^{\prime}\right]
+\left[  X,D_{11}X^{\prime}\right]  +\left[  X,D_{21}X^{\prime}\right]
\end{align*}
From the brackets of $\frak{g}_{1}\underline{\times}\frak{g}_{2}$ it follows
that
\begin{align*}
D_{11}\left[  X,X^{\prime}\right]   &  =\left[  D_{11}X,X^{\prime}\right]
+\left[  X,D_{11}X^{\prime}\right]  \\
D_{21}\left[  X,X^{\prime}\right]   &  =0\\
D_{31}\left[  X,X^{\prime}\right]   &  =\left[  D_{21}X,X^{\prime}\right]
+\left[  X,D_{21}X^{\prime}\right]
\end{align*}
This shows that $D_{11}$ is a derivation of $\frak{g}_{1}$ and that
$D_{21}\left(  C^{1}\frak{g}_{1}\right)  =0$. From the third identity it
follows moreover that $D_{31}\left(  C^{2}\frak{g}_{1}\right)  =0$.
\newline Reasoning with elements $X,X^{\prime}\in\frak{g}_{2}$ we obtain that
$D_{22}\in Der\left(  \frak{g}_{2}\right)  $ and $D_{12}\left(  C^{1}%
\frak{g}_{2}\right)  =D_{32}\left(  C^{2}\frak{g}_{2}\right)  =0$.\newline Now
let $X\in\frak{g}_{1}$ and $X^{\prime}\in\frak{g}_{2}$. Then $\left[
X,X^{\prime}\right]  \in\frak{g}_{3}$ and
\begin{align*}
D\left[  X,X^{\prime}\right]   &  =D_{13}\left[  X,X^{\prime}\right]
+D_{23}\left[  X,X^{\prime}\right]  +D_{33}\left[  X,X^{\prime}\right]  \\
&  =\left[  D_{11}X,X^{\prime}\right]  +\left[  D_{21}X,X^{\prime}\right]
+\left[  X,D_{12}X^{\prime}\right]  +\left[  X,D_{22}X^{\prime}\right]
\end{align*}
and reordering as before%
\begin{align*}
D_{33}\left[  X,X^{\prime}\right]   &  =\left[  D_{11}X,X^{\prime}\right]
+\left[  X,D_{22}X^{\prime}\right]  \\
D_{13}\left[  X,X^{\prime}\right]   &  =\left[  X,D_{12}X^{\prime}\right]  \\
D_{23}\left[  X,X^{\prime}\right]   &  =\left[  D_{21}X,X^{\prime}\right]
\end{align*}
which shows that $D_{i3}\left(  \frak{g}_{3}\right)  \subset C^{1}\left(
\frak{g}_{i}\right)  $ for $i=1,2$. \newline If $X\in\frak{g}_{1}$ and
$X^{\prime}\in\frak{g}_{3},$
\[
D\left[  X,X^{\prime}\right]  =0=\left[  X,D_{13}X^{\prime}\right]  +\left[
X,D_{23}X^{\prime}\right]
\]
and as $\left[  X,D_{13}X^{\prime}\right]  \in C^{1}\frak{g}_{1}$ and $\left[
X,D_{23}X^{\prime}\right]  \in\frak{g}_{3}$, we have
\[
\left[  X,D_{13}X^{\prime}\right]  =\left[  X,D_{23}X^{\prime}\right]  =0
\]
This proves that $D_{13}\left(  \frak{g}_{3}\right)  \subset Z\left(
\frak{g}_{1}\right)  $. Similarly we obtain $D_{23}\left(  \frak{g}%
_{3}\right)  \subset Z\left(  \frak{g}_{2}\right)  $.\newline Now, taking
$X,X^{\prime}$ generators of $\frak{g}_{1}$ and $\frak{g}_{2}$, respectively,
we have%
\[
Z\left(  \frak{g}_{1}\right)  \ni D_{13}\left[  X,X^{\prime}\right]  =\left[
X,D_{12}X^{\prime}\right]
\]%
\[
Z\left(  \frak{g}_{2}\right)  \ni D_{23}\left[  X,X^{\prime}\right]  =\left[
D_{21}X,X^{\prime}\right]
\]
This shows that $D_{12}X^{\prime}\in T_{1}$ and $D_{21}X\in T_{2}$, from which
the last assertion follows.
\end{proof}

\begin{remark}
Observe that if $T_{i}\subset C^{1}\frak{g}_{i}$, it clearly follows that
$D_{21}\left(  \frak{g}_{1}\right)  \subset C^{1}\frak{g}_{i},\newline D_{12}\left(
\frak{g}_{2}\right)  \subset C^{1}\frak{g}_{2}$ and $D_{31}\left(
C^{1}\frak{g}_{1}\right)  =D_{32}\left(  C^{1}\frak{g}_{2}\right)  =0$.
In particular, algebras satisfying this were used in [7] to prove that the deletion of the center by any derivation is not sufficient for an algebra to be characteristically nilpotent.
For an example of a CNLA $\frak{g}$ \ for which $\left\{  Y\in\frak{g\;}%
|\;\left[  Y,\frak{g}\right]  \subset Z\left(  \frak{g}\right)  \right\}
\nsubseteq C^{1}\frak{g}$ consider the algebra introduced by Luks in [10]. 
\end{remark}

The following result gives sufficient conditions to ensure that $D_{ij}\left(
\frak{g}_{j}\right)  \subset C^{1}\frak{g}_{i}\newline \left(  1\leq i,j\leq2,\;i\neq
j\right)  $ :

\begin{proposition}
Let $1\leq i,j\leq2$,\ $i\neq j$. Suppose that any generator $X$ of
$\frak{g}_{i}$ satisfies one of the following conditions 

\begin{enumerate}
\item  there exists a generator $X^{\prime}$ such that
\[
\left[  X,Y\right]  \in C^{p}\frak{g}_{i},\;\;p\geq2
\]

\item  there exist generators $X^{\prime},Y,Y^{\prime}$ of $\frak{g}_{i}$ such
that%
\[
\left[  X,X^{\prime}\right]  =a\left[  Y,Y^{\prime}\right]  \;\left(
\operatorname{mod}\;C^{2}\frak{g}_{i}\right)  ,\;\;a\neq0
\]
\end{enumerate}
\end{proposition}

Then $D_{ji}\left(  \frak{g}_{i}\right)  \subset C^{1}\frak{g}_{j}\;$and
$D_{3i}\left(  C^{1}\frak{g}_{i}\right)  =D_{3j}\left(  C^{1}\frak{g}%
_{j}\right)  =0$.

\begin{proof}
Let $X$ be a generator of $\frak{g}_{1}$ satisfying the first condition. From
the preceding proof we have $D_{31}\left[  X,X^{\prime}\right]  =\left[
D_{21}X,X^{\prime}\right]  +\left[  X,D_{21}X^{\prime}\right]  $ for
$X,X^{\prime}\in\frak{g}_{1}$.Let $X^{\prime}$ be the generator such that
$\left[  X,X^{\prime}\right]  \in C^{p}\frak{g}_{1}$ with $p\geq2$. Then
$D_{31}\left[  X,X^{\prime}\right]  =0$. If $Y,Y^{\prime}$ are arbitrary
generators of $\frak{g}_{2}$, then by the definition of $\frak{g}%
_{1}\underline{\times}\frak{g}_{2}$ it follows that $\left[  Y,X\right]
=\left[  Y^{\prime},X^{\prime}\right]  $ if and only if $X=X^{\prime}$ and
$Y=Y^{\prime}$. This shows that
\[
\left[  D_{21}X,X^{\prime}\right]  =\left[  X,D_{21}X^{\prime}\right]  =0
\]
thus $D_{21}X\in C^{1}\frak{g}_{1}$ and $D_{21}X^{\prime}\in C^{1}\frak{g}%
_{1}$.\newline The second assertion follows in a similar way.
\end{proof}

\begin{definition}
A nilpotent Lie algebra $\frak{g}_{1}$ is called an $S$-algebra if for any Lie algebra $\frak{g}_{2}$ and any derivation $D$ of $\frak{g}_{1}\underline{\times}\frak{g}_{2}$ we have 
\[
D_{21}\left(\frak{g}_{1}\right)\subset C^{1}\frak{g}_{2}
\]
\end{definition}

\begin{remark}
Clearly any CNLA $\frak{g}_{i}$ satisfying $T_{i}\subset C^{1}\frak{g}_{i}$ or the requirements of proposition $4$ is an $S$-algebra. Thus all nonsplit classical examples of CNLAs satisfy the preceding definition ( [1], [3], [9], [10]).
\end{remark}  

In what follows, we will require that the Lie algebras $\frak{g}_{i}$ are also $S$-algebras. With these conditions, we have, for any derivation $D=\sum_{i,j}D_{ij}$ :
\begin{align*}
D_{12}\left(  \frak{g}_{2}\right)   &  \subset C^{1}\frak{g}_{1}\;\text{and
\ }D_{12}\left(  C^{1}\frak{g}_{2}\right)  =0;\\
D_{21}\left(  \frak{g}_{1}\right)   &  \subset C^{1}\frak{g}_{2}\;\text{and
\ }D_{21}\left(  C^{1}\frak{g}_{1}\right)  =0\\
D_{i3}\left(  \frak{g}_{3}\right)   &  \subset Z\left(  \frak{g}_{i}\right)
,\;i=1,2\\
D_{3j}\left(  \frak{g}_{j}\right)   &  \subset\frak{g}_{3}\;\text{and
\ }D_{3j}\left(  C^{1}\frak{g}_{j}\right)  =0,\;j=1,2
\end{align*}
From this, it is easy to deduce the following relations
\begin{align}
D_{ij}D_{jk}  &  =0,\;j=1,2,\;i\neq k\\
D_{i3}D_{31}\left(  \frak{g}_{1}\underline{\times}\frak{g}_{2}\right)   &
\subset Z\left(  \frak{g}_{i}\right)  ,\;i\neq k\\
D_{ki}D_{i3}D_{3j}  &  =0,\;i,j\neq3,\;k\neq i\\
D_{i3}D_{3j}D_{jk}  &  =0,\;i,i\neq3,\;j\neq k
\end{align}

\begin{lemma}
Let $j=1,2$. For any $m\geq1$ we have
\[
D_{ij}D_{jj}^{m}D_{jk}=0,\;i\neq j,\;j\neq k
\]
\end{lemma}

\begin{proof}
As for $j\neq3$ we have $D_{jk}\left(  \frak{g}_{k}\right)  \subset
C^{1}\frak{g}_{j}$ and $D_{jj}^{m}\left(  C^{1}\frak{g}_{j}\right)  \subset
C^{1}\frak{g}_{j}$, the result follows from the fact that $D_{ij}\left(
C^{1}\frak{g}_{j}\right)  =0$.
\end{proof}

\begin{remark}
Observe that this result is false if it is not assured that $D_{ij}\left(
\frak{g}_{j}\right)  \subset C^{1}\frak{g}_{j}$ $\left(  i=1,2\right)  $.
\end{remark}

\begin{lemma}
Let $i,j\in\left\{  1,2\right\}  $. Then following equalities hold

\begin{enumerate}
\item $D_{i3}D_{33}D_{3j}\left(  \frak{g}_{1}\underline{\times}\frak{g}%
_{2}\right)  \subset Z\left(  \frak{g}_{i}\right)  $

\item $D_{i3}D_{33}D_{3j}D_{jk}=0,\;j\neq k$

\item $D_{li}D_{i3}D_{33}D_{3j}=0,\;l\neq i$
\end{enumerate}
\end{lemma}

The proof follows immediately from succesice application of the relations
$\left(  1\right)  $ to $\left(  4\right)  $. \newline From these facts, it is
elementary to see that for any $k\geq1$ we have
\begin{equation}
D^{k}=\sum_{\substack{i\neq j\\p+q=k-1}}D_{ii}^{p}D_{ij}D_{jj}^{q}%
+\sum_{\substack{i,j\neq3\\r+s+t=k-2}}D_{ii}^{r}D_{i3}D_{33}^{s}%
D_{3j}D_{jj}^{t}%
\end{equation}

\begin{remark}
Observe the similarity with the corresponding formula for the direct sum [8].
The residue of the second summand is due to the fact that $\frak{g}_{3}$ is an
abelian ideal generated by $\frak{g}_{1}$ and $\frak{g}_{2}$ and that $D_{33}$
is not a derivation.
\end{remark}

\begin{proposition}
If $D_{11}$ and $D_{22}$ are nilpotent derivations, then $D_{33}$ is a
nilpotent endomorphism.
\end{proposition}

\begin{proof}
Let $X\in\frak{g}_{1}$ and $X^{\prime}\in\frak{g}_{2}$. By the definition of
$\frak{g}_{1}\underline{\times}\frak{g}_{2}$, $\left[  X,X^{\prime}\right]
\in\frak{g}_{3}$. Then
\begin{align*}
D^{k}\left[  X,X^{\prime}\right]   &  =\left[  D^{k}X,X^{\prime}\right]
+\left[  X,D^{k}X^{\prime}\right]  \\
&  =\left[  p_{1}D^{k}X+p_{2}D^{k}X,X^{\prime}\right]  +\left[  X,p_{1}%
D^{k}X^{\prime}+p_{2}D^{k}X^{\prime}\right]
\end{align*}
It follows that
\begin{align*}
p_{1}D^{k}\left[  X,X^{\prime}\right]   &  =\left[  X,p_{1}D^{k}X^{\prime
}\right]  \;;\;p_{2}D^{k}\left[  X,X^{\prime}\right]  =\left[  p_{2}%
D^{k}X,X^{\prime}\right]  \\
p_{3}D^{k}\left[  X,X^{\prime}\right]   &  =\left[  p_{1}D^{k}X,X^{\prime
}\right]  +\left[  X,p_{2}D^{k}X^{\prime}\right]
\end{align*}
From the relations $\left(  1\right)  $ to $\left(  4\right)  $ and the form
of $D^{k}$ in $\left(  5\right)  $ it is easy to see that the last identity
simplifies to
\[
D_{33}^{k}\left[  X,X^{\prime}\right]  =\left[  D_{11}^{k}X,X^{\prime}\right]
+\left[  X,D_{22}^{k}X^{\prime}\right]
\]
Thus, if $D_{11}^{r_{1}}=D_{22}^{r_{2}}=0$, then $D_{33}^{\max\left\{
r_{1},r_{2}\right\}  }=0$.
\end{proof}

\begin{theorem}
Let $\frak{g}_{1}$ and $\frak{g}_{2}$ be two characteristically nilpotent
Lie ($S$-)algebras. Then $\frak{g}_{1}\underline{\times}\frak{g}_{2}$ is also \ a CNLA.
\end{theorem}

\begin{proof}
Let $\alpha_{1},\alpha_{2}$ be the smallest integers such that $\frak{g}%
_{i}^{\left[  \alpha_{i}\right]  }=0$. Without loss of generality suppose that
$\alpha_{2}=\max\left\{  \alpha_{1},\alpha_{2}\right\}  $. Then for any
$m\geq\alpha_{1}+2\alpha_{2}+2$ we have $D^{m}=0$. Thus any derivation of
$\frak{g}_{1}\underline{\times}\frak{g}_{2}$ is nilpotent.
\end{proof}

\begin{corollary}
Let $\frak{g}_{1},..,\frak{g}_{n}$ be characteristically nilpotent $S$-algebras. Then $\frak{g}_{1}\underline
{\times}\frak{g}_{2}\underline{\times}...\underline{\times}\frak{g}_{n}$ is a CNLA.
\end{corollary}

\begin{example}
Consider the algebra \ $\frak{g}$ of Dixmier and Lister [4] : it has basis
$\left\{  X_{1},..,X_{8}\right\}  $ and law%
\[%
\begin{array}
[c]{ccc}%
\left[  X_{1},X_{2}\right]  =X_{5} & \left[  X_{2},X_{3}\right]  =X_{8} &
\left[  X_{3},X_{4}\right]  =-X_{5}\\
\left[  X_{1},X_{3}\right]  =X_{6} & \left[  X_{2},X_{4}\right]  =X_{6} &
\left[  X_{3},X_{5}\right]  =-X_{7}\\
\left[  X_{1},X_{4}\right]  =X_{7} & \left[  X_{2},X_{6}\right]  =-X_{7} &
\left[  X_{4},X_{6}\right]  =-X_{8}\\
\left[  X_{1},X_{5}\right]  =-X_{8} &  &
\end{array}
\]
Consider the product $\frak{g}\underline{\times}\frak{g}$. It is
32-dimensional with 20-dimensional center. As this algebra satisfies 
$T=\left\{  Y\in\frak{g}\;|\;\left[  Y,\frak{g}\right]  \subset
Z\left(  \frak{g}\right)  \right\}$, its product is characteristically nilpotent. Observe in particular that $D_{33}=$. This follows from the fact that any derivation of $\frak{g}$ maps the algebra into $C^{1}\frak{g}$. Thus
\[
D\left(  \frak{g}\underline{\times}\frak{g}\right)  \subset C^{1}\left(
\frak{g}\underline{\times}\frak{g}\right)
\]
for any derivation $D$.
\end{example}

\begin{corollary}
Let $\frak{g}_{i}$ ( $i=1,2$ ) be S-algebras such that $f_{i}\left(
\frak{g}_{i}\right)  \subset C^{1}\left(  \frak{g}_{i}\right)  $ for any
$f_{i}\in Der\left(  \frak{g}_{i}\right)  $. Then $\frak{g}_{1}\underline
{\times}\frak{g}_{2}$ is characteristically nilpotent and $D\left(
\frak{g}_{1}\underline{\times}\frak{g}_{2}\right)  \subset C^{1}\left(
\frak{g}_{1}\underline{\times}\frak{g}_{2}\right)  $ for any derivation of the product.
\end{corollary}

\begin{example}
Consider the algebra $\frak{L}$ \ given by Luks in [10]. It has basis $\left\{
X_{1},..,X_{16}\right\}  $ and law%
\begin{align*}
\left[  X_{1},X_{i}\right]   &  =X_{5+i},\;i=2,3,4,5;\;\left[  X_{1}%
,X_{6}\right]  =X_{13};\;\left[  X_{1},X_{i}\right]  =X_{8+i},\;i=7,8\\
\left[  X_{2},X_{i}\right]   &  =X_{i+8},\;i=3,4,6;\;\left[  X_{2}%
,X_{5}\right]  =X_{15};\;\left[  X_{2},X_{7}\right]  =-X_{16}\\
\left[  X_{3},X_{4}\right]   &  =-X_{13}-\frac{9}{5}X_{15};\;\left[
X_{3},X_{5}\right]  =-X_{14};\;\left[  X_{3},X_{6}\right]  =-X_{16};\;\left[
X_{4},X_{5}\right]  =2X_{16}%
\end{align*}
Clearly it is an $S$-algebra, as its generators satisfy proposition $4$.
Take the product $\frak{L}\underline{\times}\frak{L}$. This algebra is
68-dimensional with 52-dimensional center of nilindex $3$. By the preceding theorem, this algebras is characteristically nilpotent.
\end{example}

\begin{problem}
It remains to see if there exist characteristically nilpotent Lie algebras $\frak{g}$ which are not $S$-algebras. Such an algebra must have at least three generators, from which at least one belongs to the transporter $T=\left\{  Y\in\frak{g}\;|\;\left[  Y,\frak{g}\right]  \subset
Z\left(  \frak{g}\right)  \right\}$ and does not satisfy the conditions of proposition $4$. 
\end{problem}

\bigskip
Let $\frak{g}$ be a characteristically nilpotent $S$-algebra with $k$
generators. Then, for any $n\geq1$, the $n^{th}$ power \underline{$\times$%
}$^{n}\frak{g}=\frak{g}\underline{\times}...\underline{\times}\frak{g}$ ( n
times )  is characteristically nilpotent. Morever, we have %

\begin{align*}
\dim\underline{\times}^{n}\frak{g}  & =\frak{}\frac{1}{2}n^{2}k^{2}+\left(
\dim\frak{g}-\frac{1}{2}k^{2}\right)  n\\
\dim Z\left(  \underline{\times}^{n}\frak{g}\right)    & =\frac{1}{2}%
n^{2}k^{2}+\left(  \dim Z\left(  \frak{g}\right)  -\frac{1}{2}k^{2}\right)  n
\end{align*}

The proof follows by induction over $n$. Observe in particular that
\[
\lim_{n\rightarrow\infty}\frac{\frac{1}{2}n^{2}k^{2}+\left(  \dim
\frak{g}-\frac{1}{2}k^{2}\right)  n}{\frac{1}{2}n^{2}k^{2}+\left(  \dim
Z\left(  \frak{g}\right)  -\frac{1}{2}k^{2}\right)  n}=1
\]
from which we obtain

\begin{theorem}
For any $\varepsilon>0$ there exists a nonsplit characteristically nilpotent
Lie algebra $\frak{h\;}$such that%
\[
\frac{codim_{\frak{h}}Z\left(  \frak{h}\right)  }{\dim Z\left(
\frak{h}\right)  }<\varepsilon
\]
\end{theorem}

where $codim_{\frak{h}}Z\left(\frak{h}\right)=dim \frak{h}-dim Z\left(\frak{h}\right)$.
\newline
Now, the number of elements
in $C^{1}\left(  \underline{\times}^{n}\frak{g}\right)  $ which are not
central equals%
\[
\dim\left(  \underline{\times}^{n}\frak{g}\right)  -\dim Z\left(  \left(
\underline{\times}^{n}\frak{g}\right)  \right)  -kn=n\;\left(  \text{co}%
\dim_{\frak{g}}Z\left(  \frak{g}\right)  -k\right)  =n\left(  \dim
C^{1}\frak{g}-\dim Z\left(  \frak{g}\right)  \right)
\]

\begin{corollary}
For any $\varepsilon>0$ there exists a nonsplit characteristically nilpotent
Lie algebra $\frak{h\;}$such that%
\[
\frac{\dim C^{1}\frak{h}-\dim Z\left(  \frak{h}\right)  }{\dim Z\left(
\frak{h}\right)  }<\varepsilon
\]
\end{corollary}

\begin{proof}
Let $N\in\mathbb{N}$ be such that $\frac{1}{N}<\varepsilon$. Let $\frak{g}$ be
any nonsplit characteristically nilpotent $S$-algebra with $k$ generators$.$
Consider the $n^{th}$-power \underline{$\times$}$^{n}\frak{g}$, where
\[
n>\frac{2N}{k^{2}}\left(  \text{co}\dim_{\frak{g}}Z\left(  \frak{g}\right)
-k\right)  -\frac{2\dim Z\left(  \frak{g}\right)  }{k^{2}}+1
\]
Then we have
\[
\frac{\dim C^{1}\left(  \underline{\times}^{n}\frak{g}\right)  -\dim Z\left(
\underline{\times}^{n}\frak{g}\right)  }{\dim Z\left(  \underline{\times}%
^{n}\frak{g}\right)  }<\frac{1}{N}%
\]
\end{proof}

\begin{remark}
Observe that for any $n,$ the number of elements in $C^{1}\left(
\underline{\times}^{n}\frak{g}\right)  $ which are not central increases
linearly, while the number of central elements increases quadratically.  
\end{remark}

\end{document}